\documentclass[12pt]{article}

\usepackage{amssymb}

\begin{document}

\newtheorem{thm}{Theorem}[section]
\newtheorem{lem}[thm]{Lemma}
\newtheorem{rem}[thm]{Remark}
\newtheorem{cor}[thm]{Corollary}
\newtheorem{prop}[thm]{Proposition}


\renewcommand{\theequation}{\arabic{section}.\arabic{equation}}
\def\proof{\noindent{\it Proof.\ }}


\newcommand{\RR}{\mathbb R}
\newcommand{\N}{{\rm I\kern-.1567em N}}
\newcommand{\Z}{{\sf Z\kern-.3567em Z}}
\def\qed{\hfill $\square$}
\newcommand{\be}{\begin{equation}} 
\newcommand{\ee}{\end{equation}}
\newcommand{\bea}{\begin{eqnarray}} 
\newcommand{\eea}{\end{eqnarray}}
\newcommand{\bean}{\begin{eqnarray*}} 
\newcommand{\eean}{\end{eqnarray*}}
\newcommand{\rf}[1]{(\ref {#1})}
\newcommand{\un}{\rm 1\!\!I}
\newcommand{\la}{\langle}
\newcommand{\ra}{\rangle}
\newcommand{\iom}[1]{\int_\Omega #1\, dx}
\newcommand{\io}[1]{\int_\Omega #1}
\newcommand{\idom}[1]{\int_{\partial\Omega} #1\; d\sigma}
\def\th{\theta}
\def\vt{\vartheta}
\def\div{\nabla\cdot}
\def\f{\varphi}
\def\r{\rightarrow}
\def\j{{\bf j}}
\def\h{{H^1(\Omega)}}
\title{On an evolution system describing\\
self-gravitating particles in microcanonical setting} 
\author{Robert Sta\'nczy\\
\small Instytut Matematyczny, Uniwersytet Wroc{\l}awski, \\
\small pl. Grunwaldzki 2/4, 50--384 Wroc{\l}aw, Poland;\\
\small{\tt stanczr@math.uni.wroc.pl}}

\date{\today}
\maketitle

\begin{abstract} 
The global in time existence of solutions of a system describing 
the interaction of gravitationally attracting particles with
a general diffusion term and fixed energy is proved.
The presented theory covers the case of the model with diffusion 
that obeys Fermi--Dirac statistics. Some of the results apply to
the dissipative polytropic case as well.
\end{abstract}





\baselineskip=16pt


\section{Introduction}

\setcounter{equation}{0}

We consider the following initial-boundary value problem
\bea
\label{ibp1}
n_t = \nabla\cdot\left(D\left(\nabla p(n) + 
  n \nabla\varphi \right)\right) & \;\;\mbox{ in }\;\; & \Omega\times
(0,\infty)\,,\\
\label{ibp2}
\Delta\varphi = n & \;\;\mbox{ in }\;\; & \Omega\times (0,\infty)\,,\\
\label{ibp3}
\left(\nabla p(n) + \nabla\varphi \right)\cdot\bar{\nu} =
\varphi = 0 & \;\;\mbox{ on }\;\; & \partial\Omega\times (0,\infty)\,,\\
\label{ibp4}
n(0) = n_0 \geq 0 & \;\;\mbox{ in }\;\; & \Omega\subset \RR^d \,,
\eea
with the pressure in the self-similar form
\be \label{ssf}
p(n,\th)=\th^{d/2+1}P(n\th^{-d/2})
\ee
for a given function $P$ and some diffusion coefficient 
$D,$ motivated by applications to statistical mechanics
and describing self--attracting clouds of particles modelling elliptical
gallaxies, globular clusters, interstellar medium or cores of neutron stars
among others (cf. \cite{CRS} and references therein).
These sorts of problems were considered among others in \cite{BS-d}
including numerous pressure formulae coming from statistical mechanics:
Maxwell--Boltzmann, Fermi--Dirac, Bose--Einstein and polytropic distributions.
The common feature of all these examples is the self-similar profile
of the pressure (\ref{ssf}). In this paper we focus our attention on
the Fermi--Dirac model although we formulate the results in a
more general setting. The pressure in this model assumes an intermediate
form between the linear Maxwell--Boltzmann case at zero and a polytropic, power--like
form at infinity. In \cite{BLN} the authors proved the local and global 
existence for the specific choice of the diffusion parameter $D$ 
corresponding to the Fermi--Dirac statistics in the isothermal case, 
i.e., with fixed, constant temperature $\theta$.
Moreover, the asymptotic behaviour with the possibility of the evolution
towards steady states was adressed therein (cf. also \cite{S-r}), 
whereas in \cite{BNS} 
some results for nonisothermal case were established. In the 
aforementioned papers, some {\emph{a~priori}} estimates for the
density $n$ and the pressure $P$ were also provided. 
For physical motivations one can see the series of papers of 
Chavanis and collaborators including \cite{Ch}-\cite{CSR}.
Related Keller-Segel model in mathematical biology 
was recently studied in \cite{LS}, \cite{Su} 
and \cite{SK} and the blow-up for large data was proved.

Note that, as a consequence of
\rf{ibp3}, total mass
\be
M=\iom{n(x,t)}\label{M}
\ee
is conserved during the evolution of the system. 

In the first part of the paper we will extend the results of \cite{BLN} 
to allow more general pressure and nonconstant temperature for 
the dimension $d\in [2,4]$ with small mass $M$. The local and global
results for the parabolic perturbation is contained in Proposition \ref{loc},
and the result for the original parabolic-elliptic problem can be found 
in Theorem \ref{thglobal}. 
As far as steady states are
concerned  nonexistence results 
hold for $d>2(1+\sqrt{2})$ (cf. \cite{S-d}) we can expect that 
global existence result holds only for the dimension $d\le 4$ 
when global minimizer for the entropy functional is attained 
as proven in \cite{S-d}, \cite{S-p} or \cite{BLN}. In fact analogous nonexistence
results also hold for a problem related to \rf{ibp1}--\rf{ibp4} but
with constant diffusion parameter $D$ (cf. \cite{BLN}). Thus 
one can conjecture that the gap $d\in (4,2(1+\sqrt{2}))$
is left for the existence of the critical points (possibly unstable)
of another type than the extremal ones.
 
Next, we shall use the aforementioned existence theorems
for a given temperature $\theta(t)$ at time $t$ to prove the existence theorem
\ref{final} in the microcanonical (nonisothermal) setting, i.e. with the given energy
and the temperature to be determined so that the energy 
relation \rf{ene} is satisfied. Steady states for the model were 
considered, among others, in \cite{S-m}.
Thus, we will show that in low dimensions $2\le d \le 4$ 
for small mass and domination of the thermal energy
a gravo-thermal catastrophe (white dwarfs in a physical interpretation) 
does not occur for this system, i.e., 
neither blow-up for the density nor the vanishing of the 
temperature takes place.

Finally, in Appendix, we gather the properties of some special 
functions appearing in the Fermi--Dirac model.

First, notice that due to the self-similar structure of the pressure
(\ref{ssf}) for the specific canonical diffusion coefficient $D=P'$
the system (\ref{ibp1})--(\ref{ibp4}) can be transformed to the following
one (cf. also the Appendix for Fermi--Dirac case and the papers 
\cite{CRS}, \cite{BLN}, \cite{BNS} where such $D$ was used). 
Thus we arrive at the system
\bea
\label{bvp1}
n_t = \nabla\cdot\left(\th P'^2\ \nabla n + nP'\
  \nabla\varphi \right) &  \;\;\mbox{ in }\;\; & \Omega\times
(0,\infty)\,,\\ 
\label{bvp2}
\Delta\varphi = n & \;\;\mbox{ in }\;\; & \Omega\times (0,\infty)\,,\\
\label{bvp3}
\left(\th P'^2\ \nabla n + nP'\
  \nabla\varphi \right)\cdot \bar{\nu} =
\varphi = 0 & \;\;\mbox{ on }\;\; & \partial\Omega\times (0,\infty)\,,\\ 
\label{bvp4}
n(0) = n_0 \ge 0 & \;\;\mbox{ in }\;\; & \Omega\subset \RR^d\,,
\eea
where we  suppose that $2\le d\le 4$, 
the temperature is a fixed continuous function
$\th:[0,\infty)\rightarrow [a,b],$ with some positive numbers $a$ and $b$
with values to be determined later.
Morover, we look for the solutions of (\ref{bvp1})--(\ref{bvp4})
satisfying the energy relation given by
\be \label{E}
E=\frac{d}{2}\iom{\th^{d/2+1}P(n\th^{-d/2})}+\frac{1}{2} \iom {n\varphi}
= {\rm \, const}.
\ee

The steady state problem with the prescribed energy for the linear diffusion was 
considered, among others, in \cite{BDEMN}.
The main results of this paper (to be specified in the next sections) 
can be stated as follows
\begin{thm}
If $P$ is the Fermi--Dirac pressure and mass $M>0$ is sufficiently small
the problem \rf{bvp1}-\rf{bvp4} admits at least one global weak solution 
for $d\le 3$ and a local one if $d=4$ for a given continuous function
$\theta(t).$
Moreover, there exists a local weak solution to the problem 
\rf{bvp1}-\rf{bvp4} with the energy given by \rf{E} for $d\le 4$.
\end{thm}

Now we shall sketch the method of proving the above theorem.

First, we regularize the problem to obtain a parabolic system 
and to apply general Amann theory. Next sign--sensitive 
a priori bounds together with a bootstrap argument are used
to prove global or local existence depending on the dimension $d$. 
Then we go with the parameter 
to infinity and obtain the corresponding existence 
result in weak sense for the original elliptic--parabolic system.

We introduce a new temperature, call it $\vt,$ defined 
implicitly by the aforementioned energy relation, i.e.
\be\label{ene}
E=\frac{d}{2}\iom{\vt^{d/2}P(n\vt^{-d/2})}+
\frac{1}{2} \iom{n\varphi} = {\rm \, const}.
\ee
Note the implicit dependence of the `new' temperature $\vt$ on 
the old one $\th$ via $n,\f$ (solving (\ref{bvp1})--(\ref{bvp4}) for given $\th$) 
in the above formula.
In sections 2 and 3, for given $\th$, 
we solve (\ref{bvp1})--(\ref{bvp4}) to get $n, \varphi.$ 
Then, in section 4 for given value of the energy 
$E,$ we use the implicit formula \rf{ene} for $\vt$ and ask whether the 
operator $\mathcal{T}:\th\mapsto\vt$ defined by \rf{ene} 
has a fixed point. The problem of {\emph{a~priori}}
bounds for the temperature $\th,$ determined by \rf{E}, was
addressed in \cite{BS-t}. In the last section the properties of the special
Fermi--Dirac pressure function have been gathered.

\medskip

\noindent {\bf Notation.} By $C$ we will denote inessential constants, which may vary from
one line to another. By $|\cdot|_p,$ for $p\geq 1,$ 
we shall denote the standard $L^p(\Omega)$ norm. 
By smoothness we shall always mean $C^2$ regularity and it will apply only to 
the function $P$ and is explicitly stated at the beginning of the next section. 
Finally, both $\Vert \cdot \Vert_{H^1}$ and $\Vert \cdot \Vert_{*}$ will denote
the norm in the Sobolev space $H^1(\Omega)$ with $L^2(\Omega)$ and $L^1(\Omega)$
term correspondingly.
\section{The existence result for the perturbation}

In this section we follow the lines of the proof of the existence
addressed in \cite{BLN}, where the authors considered a specific Fermi--Dirac 
density $P=P_{FD}$ defined by \rf{PFD}, dimension $d=3$ and a constant temperature 
$\th,$ whereas here we will just exploit 
smoothness of the pressure $P\in C^2\left([-\delta,\infty);[0,\infty)\right),$ 
the crucial estimates: $a\le \th(t)\le b$ with some $a>0,b>0$, and for 
$z\ge 
0$,
\bea
\max\{p_0,p_1z^{2/d}\}\le P'(z) \le p_2(1+z^{2/d}), \label{p0} \\
zP''(z) \le p_3 (1+z^{2/d}). \label{p00}
\eea
These assumptions imply that, changing $p_2$ if necessary,
\bea
\max\{p_0z,p_1z^{1+2/d}\} \le P(z) \le p_2(1+z^{1+2/d}), \label{p1} \\
zP'(z)\le CP(z), \label{p2} \\
zP''(z)\le CP'(z). \label{p3}
\eea

\medskip

In order to study the well-posedness of (\ref{bvp1})--(\ref{bvp4}),
for $k\ge 1$ and $P'=P'(n\th^{-d/2}),$
we consider the following regularized initial-boundary value problem
\bea
\label{reg1}
n_t = \nabla\cdot\left(\th P'^2\ \nabla n + nP'\
  \nabla\varphi \right) & \;\;\mbox{ in }\;\; & 
\Omega\times
(0,\infty)\,,\\ 
\label{reg2}
\varphi_t - k\ \Delta\varphi = -k\ n & \;\;\mbox{ in }\;\; &
\Omega\times (0,\infty)\,,\\ 
\label{reg3}
\left(\th P'^2\ \nabla n + nP'\
  \nabla\varphi \right)\cdot\bar{\nu} =
\varphi = 0 & \;\;\mbox{ on }\;\; & \partial\Omega\times
(0,\infty)\,,\\ 
\label{reg4}
(n(0),\varphi(0)) = (n_0,\varphi_0) & \;\;\mbox{ in }\;\; &
\Omega\subset \RR^d\,,
\eea
with
\be
\label{ini}
n_0\in \mathcal{C}^\infty(\bar{\Omega}) \;\;\mbox{
  such that }\;\; n_0\ge 0\,, \qquad M=\vert n_0\vert_1, \;\;\mbox{ and
  }\;\; \varphi_0=0\,.
\ee

\medskip

For this parabolic system we first use the theory developed by Amann \cite{Ama} to prove the
local well-posedness of (\ref{reg1})--(\ref{reg4}) and then the global one. The proposition
formulated below and its proof has been adapted from \cite{BLN} to cover the case of 
variable temperature $\th(t),$ slightly more general pressure $P$ than the 
Fermi--Dirac one, and any dimension $2\le d \le 4$.

\begin{prop}\label{loc}
If $d\ge 2$ and the function $P$ is smooth and satisfies $P'\ge p_0>0$ 
and the temperature is continuous and satisfies $\th(t)\ge a$ then the initial-boundary value 
problem 
(\ref{reg1})--(\ref{reg4}) \- has a unique maximal classical solution  
$$
(n,\varphi) \in \mathcal{C}(\bar{\Omega}\times [0,T_{max});\RR^2)\cap
\mathcal{C}^{2,1}(\bar{\Omega}\times (0,T_{max});\RR^2) 
$$ 
for some $T_{max}\in (0,\infty]$. In addition,
$n(t)\ge 0 \;\mbox{ for }\; t\in [0,T_{max})\,$. 
\medskip

Furthermore, $T_{max}=\infty$ if there are $\varepsilon>0$ and a locally
bounded function $\omega: [0,\infty)\to [0,\infty)$ such that, for
every $T>0$, the estimate holds  
\be
\label{global}
\Vert n(t)\Vert_{\mathcal{C}^\varepsilon} +
\Vert\varphi(t)\Vert_{\mathcal{C}^\varepsilon}\le \omega(T) \;\;\mbox{
  for }\;\; t\in [0,T_{max})\cap [0,T]\,.  
\ee
It is the case for $d\le 3$ if $P$ additionally satisfies \rf{p0}--\rf{p00} and  
\rf{PRp} with $\vert R'(z)\vert z^{1/2-1/d}\leq B$ and $\theta(t) \le b$
for some constants $b>0,B>0.$
For $d=4$ we claim only the local existence result since the boostrap argument
does not yield \rf{global}.
\end{prop}

\noindent\textbf{Proof.} 
We set $D_0=(-\delta,\infty)\times\RR$, $u=(n,\varphi)$ with
$u_0=(n_0,\varphi_0)$, and by the assumptions, define 
$a\in\mathcal{C}^2(D_0;\mathcal{M}_2(\RR))$ and $f\in
\mathcal{C}^2(D_0;\RR^2)$ by  
$$ a(u) = \left(
\begin{array}{cc}
\th\left(P'(n\th^{-d/2})\right)^2 & nP'(n\th^{-d/2}) \\
0 & k 
\end{array}
\right)\,, \qquad f(u) = \left(
\begin{array}{c}
0 \\
- k\ n 
\end{array}
\right)\,.$$
Next, for $v\in D_0$, we introduce the operators
\bean
\mathcal{A}(v)u & = & - \sum_{i=1}^d \sum_{j=1}^d \partial_i \left(
  a_{ij}(v)\ \partial_j u \right)\,, \\ 
\mathcal{B}(v)u & = &  b\ \sum_{i=1}^d \sum_{j=1}^d \bar{\nu}_i
\left( a_{ij}(v)\ \partial_j u \right) + (I_2-b)\ u\,,  
\eean
where $a_{ij}(v) = a(v)\ \delta_{ij}$, $1\le i, j\le d$, and 
$$
I_2= \left( 
\begin{array}{cc}
1 & 0 \\
0 & 1
\end{array} 
\right) \,, 
\qquad b= \left( 
\begin{array}{cc}
1 & 0 \\
0 & 0
\end{array} 
\right) \,.
$$
Then, an abstract formulation of (\ref{reg1})--(\ref{reg4}) reads 
\bea
\label{abs1}
u_t + \mathcal{A}(u)u & = & f(u)\,,\\
\label{abs2}
\mathcal{B}(u)u & = & 0\,,\\
\label{abs3}
u(0) & = & u_0\,.
\eea
Thanks to the strict positivity of $P'(z)\ge p_0$ and the lower bound for
the temperature $\th(t)\ge a$, the eigenvalues of the matrix $a(v)$
are positive for each $v\in D_0$, and the boundary-value operator
$(\mathcal{A},\mathcal{B})$ is of separated divergence form 
and is normally elliptic in the sense of \cite[Section 4]{Ama}
Therefore we may apply \cite[Theorem~14.4 and
Theorem~14.6]{Ama} to conclude that, 
for some $T_{max}\in (0,\infty],$ (\ref{abs1})--(\ref{abs3})
has a unique maximal classical solution
$$ u=(n,\varphi) \in \mathcal{C}(\bar{\Omega}\times [0,T_{max});D_0)
\cap \mathcal{C}^{2,1}(\bar{\Omega}\times (0,T_{max});D_0).$$ 
Also, since $n_0\ge 0$ and the first 
component of $f(u)$ is equal to zero, the comparison principle (see, e.g.,
\cite[Theorem 15.1]{Ama} or \cite[Corollary~I.2.1]{LSU}) implies that 
$n(t)\ge 0$ for $t\in [0,T_{max})$. Furthermore, since $f$ does not depend
on $\nabla u$ and $n\ge 0$, Theorem~15.3 in \cite{Ama} ensures that
$T_{max}=\infty$  if there are $\varepsilon>0$ and a locally bounded
function $\omega:  [0,\infty)\to [0,\infty)$ such that (\ref{global})
holds true for every $T>0$. The fact that the assumptions imposed on $P$ 
guarantee \rf{global} requires some preparatory lemmas and is postponed until the end 
of this section. \hfill $\square$

\medskip

We proceed to present a series of lemmas which will guarantee 
that (\ref{global}) is satisfied for $2\le d\le 3$ and thus ascertain
the global solvability of the perturbed problem. 
We recall after \cite{BS-d} that the neg-entropy functional $\mathcal{W}$ 
$$
\mathcal{W} = \int_\Omega \left( nH -
  \left( \frac{d}{2}+1\right)\ P\th^{d/2} \right)\ dx\, 
$$
plays the role of a Lyapunov functional for the original and regularized problem. 
The function $H(z)$ depending on $z=n\th^{-d/2}$ is a primitive of
$P'(z)/z.$ However, in our case this functional is not useful for {\emph{a~priori}} 
estimate of the density $n$ (contrary to isothermal case \cite{BLN}) 
as due to (\ref{PRp}) it is of too low order in $n$. 
Indeed, the order is $1-2/d$ for the Fermi--Dirac case
to be exact (cf. \cite[Lemma 3.6]{BNS}), 
and as such does not provide any reasonable 
{\emph{a~priori}} estimates for the density $n$. 
On the other hand, it can be used to get {\emph{a~priori}} 
bounds for the fixed points 
of the temperature operator ${\mathcal T}$ as was done in \cite{BS-d} and is 
presented in section 4.
In the isothermal ($\th=const$) case a crucial $L^{1+2/d}$ bound was obtained
from the fact that the entropy (other than ${\mathcal W}$) was coercive in this space.
As one can see it is not the case for ${\mathcal W}$. For the details of the
nontrivial derivation of the entropy ${\mathcal W}$ one can see \cite{BS-d} 
and for its application to get {\emph{a~priori}} bounds for the temperature - 
\cite{BS-t} 
and \cite{S-b}.

Now, we are going to formulate analogous results to the ones presented in \cite{S-b}
where {\emph{a~priori}} bounds for the limit parabolic-elliptic system, as 
$k\r\infty$, were obtained. 

\begin{lem}\label{le3p}
Assume that, for $d\ge 2$,
\be \label{PRp}
P(z)=p_1z^{1+2/d}+R(z),
\ee
where the lower order term satisfies $\vert R'(z)\vert z^{1/2-1/d}\leq B.$
Then, for any fixed $T>0$ and any $t\in [0,T]\cap [0,T_{max})$, the following growth 
condition holds
\be \label{e}
\frac{d}{dt}\left( \frac{d}{2} \io{p_1n^{1+2/d}} + \io{n\f } 
+ \frac{1}{2} \io{\vert \nabla \varphi\vert^2} \right) 
+ \frac{1}{k}\io{\f_t^2} \leq 
C\theta^{d/2}  \io{\vert \nabla \varphi\vert^2}.
\ee
\end{lem}
\noindent\textbf{Proof.} Let $t\in [0,T]\cap [0,T_{max})$ and recall
that both $P$ and $P'$ are the functions of $n\th^{-d/2}.$
Now, we multiply (\ref{reg1}) by $\frac{d}{2}n^{2/d}$ and
integrate over $\Omega$ to obtain
$$
\frac{d^2}{2(d+2)}\frac{d}{dt} \iom{n^{1+2/d}} =
-  \th \iom{P'^2\vert\nabla n\vert^2n^{2/d-1}}
- \iom{ P'n^{2/d}\nabla n\cdot \nabla\varphi}\,.
$$
Similarly, multiplying (\ref{reg1}) by $A\varphi,$ we get
$$
A\iom{n_t\f }=
- A \iom{P'\ n\ \vert\nabla \varphi\vert^2}
 - A \th \iom{P'^2\ \nabla n\cdot \nabla\varphi}\,.
$$
Summing up the above equalities and using the H\"older
inequality
\bean
&&\left\vert \iom{P'\th^{1/2}\ n^{1/d-1/2}\nabla n \cdot \
\nabla\varphi
\left( \th^{-1/2}n^{1/d+1/2}+AP' \th^{1/2}n^{-1/d+1/2}\right)}
\right\vert\\
&&\leq \iom{ P'^2\th n^{2/d-1}\vert \nabla n\vert^2+
\frac{1}{4}\vert \nabla
\varphi\vert^2\left(\th^{-1/2}n^{1/d+1/2}
+AP'\th^{1/2}n^{-1/d+1/2}\right)^2}\,,
\eean
taking $A(d+2)p_1=d$ and 
$$
\iom{n_t\f }=\frac{d}{dt}\left( \iom{n\f +\frac{1}{2}\vert \nabla \varphi\vert^2} \right) 
+ \frac{1}{k}\iom{\varphi_t^2},
$$
we arrive at
\bean
\frac{d}{dt} \left( \frac{d}{2} \iom{p_1\,n^{1+2/d}} + \iom{n\f } 
+ \frac{1}{2} \iom{\vert \nabla \varphi\vert^2} \right) 
+ \frac{1}{k}\iom{\f_t^2} \\
\leq \frac{1}{4A}\iom{\vert \nabla \varphi\vert^2 
\left(AP'\th^{1/2}n^{-1/d+1/2}-\th^{-1/2}n^{1/d+1/2}\right)^2}\,.
\eean
This yields the claim, with $C=\frac{dB^2}{4(d+2)p_1}$, from the
assumption on $R'$ applied to
the differentiated pressure $P'(z)=(p_1(d+2)/d) z^{2/d}+R'(z).$ \qed

\medskip

\noindent {\bf Remark.} Note that the above theorem holds both
in the polytropic case with $R(z)=0$ and, less obviously, in the Fermi--Dirac case
as explained below. Indeed, by the properties of Fermi functions (cf. Lemmma 
\ref{FAS} from the Appendix or for more properties see \cite[Sec.5]{BNS}) 
we get $\vert R'(z)\vert \leq Bz^{-2/d}$ at $z=\infty$ and $\vert R'(z)\vert \leq B$ 
at $z=0$, implying the required estimate if $2\le d \le 6.$

\medskip

Next lemmas will allow us to estimate the right hand side of (\ref{e}).
\begin{lem}\label{Pbd}
For any $2 \leq d \leq 4$ we have 
\begin{equation}\label{pfi}
\left|\iom{n\f}\right|\le
CM^{1/2-1/d}\left(\iom{n^{1+2/d}}+ \iom{\vert \nabla \varphi\vert^2}\right).
\end{equation}
\end{lem}
\noindent\textbf{Proof.}
The proof of \rf{pfi} involves standard H\"older and Sobolev--Gagliardo--Nirenberg
inequalities as follows
\bean
&&\left|\iom{n^{1/2-1/d}n^{1/2+1/d}\f}\right|\le 
M^{1/2-1/d} \left( \iom{n\f^{\frac{2d}{d+2}}} \right)^{\frac{d+2}{2d}} \\
&&\le M^{1/2-1/d} |n|_{\frac{d+2}{4}}^{\frac{d+2}{2d}} |\f|_{\frac{2d}{d-2}} \leq 
CM^{1/2-1/d}\left( \iom{n^{1+2/d}}+ \iom{\vert \nabla \varphi\vert^2}\right)
\eean
due to the inequality $\frac{d+2}{4} \leq \frac{d+2}{d}$ and the fact 
that $H^1_0(\Omega)$ can be 
imbedded in $L^{\frac{2d}{d-2}}(\Omega).$ The proof of the case $d=2$ is 
straightforward by the Poincar\'e inequality. \qed

\medskip

In low dimensions a similar argument leads to another estimates (cf. \cite{S-b}).

\begin{lem}\label{Pbd3}
For $d=2$ we have

\begin{equation}\label{pfi2}
\left|\iom{n\f}\right|\le
C\iom{n^2}\,,
\end{equation}
while for $d=3$ the estimate
\begin{equation}\label{pfi3}
\left|\iom{n\f}\right|\le
CM^{7/3}+\frac{d}{2}\iom{n^{5/3}}\,,
\end{equation}
holds.
\end{lem}

\medskip

Now, we are ready to deduce the following lemma on {\emph{a~priori}} estimates.
\begin{lem}\label{le4p}
Assume that $2\leq d \leq 4,$ condition \rf{PRp} holds for a smooth function $P$
satisfying $|R'(z)|\le Bz^{1/d-1/2}$ 
and the temperature is bounded from above $\theta(t)\le b.$ 
Then for any $t\in [0,T]\cap [0, T_{max})$
and sufficiently small data, i.e. mass $M$ if $2<d\le 4$ or the Poincar\'e constant
for $d=2$, we have
\be \label{es3p}
\frac{d}{2} \iom{p_1\,n^{1+2/d}} + \iom{n\f } + \frac{1}{2} \iom{\vert \nabla \varphi\vert^2} 
+ \frac{1}{k}\int_0^t\iom{\f_t^2}ds \leq C.
\ee
Moreover, each of the terms appearing on the left hand side of the above inequality is bounded
and the constant $C$ may depend on the initial data. If $d=3$ the assumption on 
the smallness of $M$ can be relaxed due to Lemma \ref{Pbd3}.
\end{lem}
\noindent\textbf{Proof.}
Starting with the direct consequence of Lemma \ref{Pbd}, true for sufficiently small mass $M$ and 
large $C$  (if $d=2$, instead of making mass $M$ small, we have to assume that the 
constant from the Poincar\'e inequality is smaller than $2$),
$$
\iom{\vert \nabla \varphi\vert^2} \leq C \left( \frac{d}{2} \iom{p_1\,n^{1+2/d}} + 
\iom{n\f } + \frac{1}{2} \iom{\vert \nabla \varphi\vert^2} \right)
$$
we plug this into (\ref{e}) and integrate with respect to time to arrive, with 
possibly a larger $C$ dependent on the initial value of the right hand 
side of the above inequality, at
$$
\frac{d}{2} \iom{p_1\,n^{1+2/d}} + \iom{n\f } + \frac{1}{2} \iom{\vert \nabla \varphi\vert^2} 
\leq C \exp{\left( \int_0^t C\th(s)^{d/2} ds\right) }
$$
which together with an upper bound on $\th$ ends the proof. 
Note that \rf{pfi} from Lemma \ref{Pbd} shows that the negative term $\iom{n\f }$ is dominated 
by the positive ones and thus the last claim of the lemma is ascertained. \qed

\medskip

The integral version of the estimate \rf{e} from Lemma \ref{le3p} follows by an argument 
similar to the one in the proof of Lemma \ref{le4p} and reads

\be \label{ii}
V(t) \le V(0) \exp \left(C \int_0^t \th(s)^{d/2} ds\right)\,,
\ee
where 
$V(t) \stackrel{\rm df} = \frac{d}{2} \io{p_1n^{1+2/d}} + \io{n\f }
+ \frac{1}{2} \io{\vert \nabla \varphi\vert^2}$. Note that for $d=3$ 
we can add $CM^{7/3}$ in the definition of $V$ and thus relax the
assumption on smallness of mass $M$.

\medskip

Now, we state similarly as in \cite{BLN}, where only three--dimensional case was 
treated, a lemma on the improved regularity of $\nabla\varphi$.

\begin{lem}\label{le2}
Let $q, \alpha\in (1,\infty)$, $d\geq 2$ and $T>0$. There is a constant $C$
depending on $q, \alpha, d$ and $T$ such that, for $t\in [0,T_{max})\cap [0,T],$
\be \label{es2}
\int_0^t \vert\nabla\varphi(s)\vert_{\alpha}^q\ ds \le C\vert n
\vert^q_{L^q(0,t;L^{\frac{d\alpha}{d+\alpha}}(\Omega))}. 
\ee
\end{lem}
\noindent\textbf{Proof.} 
We infer from \cite[Corollaire~1.1]{Lam}, as in \cite{BLN} where for 
$\alpha=\frac{d(d+2)}{d^2-d-2}$ and $d=3$ the authors used the bound with 
$d\alpha/(d+\alpha)=1+2/d$ norm of $n$, that
$$
\frac{1}{k}\ \vert \varphi_t
\vert_{L^q(0,t;L^{\frac{d\alpha}{d+\alpha}}(\Omega))} + \vert \Delta\varphi
\vert_{L^q(0,t;L^{\frac{d\alpha}{d+\alpha}}(\Omega))} \le C\ \vert n
\vert_{L^q(0,t;L^{\frac{d\alpha}{d+\alpha}}(\Omega))}\,. 
$$
Now, if $t\in [0,T_{max})\cap [0,T]$, we get from the above inequality
$$
\int_0^t \Vert\varphi(s)\Vert_{W^{2,\frac{d\alpha}{d+\alpha}}(\Omega)}^q \ 
ds \le C\ \vert n \vert^q_{L^q(0,t;L^{\frac{d\alpha}{d+\alpha}}(\Omega))}\,.
$$
To conclude we use the imbedding of $W^{2,\frac{d\alpha}{d+\alpha}}(\Omega)$
in $W^{1,\alpha}(\Omega)$. 
\qed

\medskip

Furthermore, an $L^2$-estimate is available for $n$.

\begin{lem} \label{le4}
Let $T>0$, $2\le d\le 4$, $P$ be smooth and satisfy \rf{p0}
and (\ref{PRp}).
There are constants $c,C>0$ depending on sufficiently 
small mass $M$, bounds on $\theta(t)\in [a,b]$ and $P$ 
and the initial data such that,
for $t\in [0,T_{max})\cap [0,T]$,
\be \label{es3}
\vert n(t)\vert_2^2 + c\int_0^t |\nabla n(s)|_{2}^2\ ds \le C
\ee
In fact, a constant $C$ is a function
of the integral $\int_0^t |n(s)|_{1+2/d}\ ds,$
locally bounded in $t$, which can be estimated by constant due
to \rf{es3p}.
\end{lem}
\noindent\textbf{Proof.} Note that the estimate of the first term 
in \rf{es3} follows from Lemma \ref{le4p} if $d=2$.
Let $t\in [0,T_{max})\cap [0,T]$ and multiply
(\ref{reg1}) by $2 n$, and integrate over $\Omega$ to obtain
$$
\frac{d}{dt} \vert n(t)\vert_2^2 + 2 \th\ \int_\Omega P'^2\
\vert\nabla n\vert^2\ dx = - 2\ \int_\Omega nP'\ \nabla n\cdot
\nabla\varphi\ dx\,.
$$
Next, we have
$$
2\ \left\vert \int_\Omega nP'\ \nabla n\cdot \nabla\varphi\ dx
\right\vert \le \th\ \int_\Omega P'^2\ \vert \nabla n\vert^2\
dx + \frac{1}{\th}\ \int_\Omega n^2\ \vert\nabla\varphi\vert^2\
dx
$$
by the Young inequality, whence
\be
\label{tt1}
\frac{d}{dt} \vert n(t)\vert_2^2 + \th\ \int_\Omega P'^2\
\vert\nabla
n\vert^2\ dx \le \frac{1}{\th}\ \int_\Omega n^2\
\vert\nabla\varphi\vert^2\ dx\,.
\ee
For $d=2$, by Lemma \ref{le2} and Lemma \ref{le4p}, we deduce that
\begin{equation}\label{fii}
|\nabla \varphi|_{\infty}\le C 
\end{equation}
and that the integrated with respect to
the time variable the right hand side of \rf{tt1} is bounded. Thus the 
estimate \rf{es3} is proved in this case. If $d\ge 3$ a longer argument is
required. Namely, it follows from the H\"older and Young inequalities that for any 
$\varepsilon>0, \alpha>2$ and some $C=C_\varepsilon$
$$
\int_\Omega n^2\ \vert\nabla\varphi\vert^2\ dx \le \vert
n^2\vert_{\alpha/(\alpha-2)}\ \vert\nabla\varphi\vert_{\alpha}^2 \leq
\frac{\varepsilon p_1^2}{3}\left\vert n^2
\right\vert_{\alpha/(\alpha-2)}^{\frac{\alpha (d+6)}{d(\alpha+2)}} + C
\vert\nabla\varphi\vert_{\alpha}^{\frac{2\alpha (d+6)}{6\alpha-2d}}\,.
$$
Then, interpolating with positive $\beta=\frac{4(2\alpha-d-2)}{d(\alpha+2)},$ we get 
\be\label{nin}
\vert n^2 \vert_{\alpha/(\alpha-2)}^{\frac{\alpha(d+6)}{d(\alpha+2)}} 
\leq M^{\beta} \vert n^{1+2/d}\vert^2_{2d/(d-2)}
\ee
and $P(z) \ge p_1 z^{1+2/d}$, or precisely $\th^{d/2+1}P(n\th^{-d/2}) 
\ge p_1 
n^{1+2/d}$, implies that
\bea
\int_\Omega n^2\ \vert\nabla\varphi\vert^2\ dx \le \frac{\varepsilon}{3}M^{\beta} \vert
\th^{d/2+1}P\vert_{2d/(d-2)}^{2}
+ C \vert\nabla\varphi\vert_{\alpha}^{\frac{2\alpha (d+6)}{6\alpha-2d}} \\
\le \frac{1}{3}M^{\beta}\Vert \th^{d/2+1} P\Vert_{*}^{2}
+ C\ \vert\nabla\varphi\vert_{\alpha}^{\frac{2\alpha (d+6)}{6\alpha-2d}}\,,
\eea
where the last inequality follows from the continuous imbedding of
$H^1(\Omega)$ in $L^{2d/(d-2)}(\Omega)$ with a constant $\varepsilon^{-1}$ and 
the norm 
\begin{equation}\label{nor}
\Vert z\Vert_{*}^2=\iom{|\nabla z|^2}+\left(\iom{|z|}\right)^2.
\end{equation}

Consequently, by (\ref{tt1}), it follows that, for 
$c=1-\frac{1}{3}M^\beta,$
\be\label{dnl}
\frac{d}{dt} \vert n(t)\vert_2^2 + c\int_\Omega
\th P'^2\vert \nabla n\vert^2\ dx \le \frac{C}{\th}\ \left( \vert \th^{d/2+1}P
\vert_1^2 + \vert\nabla\varphi\vert_{\alpha}^{\frac{2\alpha 
(d+6)}{6\alpha-2d}} \right)\,.
\ee
Now, integration with respect to  time, the assumption on growth of $P$ 
and Lemma \ref{le2} with $q={\frac{2\alpha (d+6)}{6\alpha-d}}, \alpha=1+d/2, \beta=0,$ 
for $2<d<(3+\sqrt{17})/2$ (including $d=3$) yields the estimate by the time
integral of $|n|_{1+2/d}$.
Finally, due to Lemma \ref{le4p} providing a bound for $|n|_{1+2/d}$, 
the estimate (\ref{es3}) is proven for any mass $M$. 

Now, allowing higher dimensions $2<d<2(1+\sqrt{2})$ (including $d=4$) 
we use $q={\frac{2\alpha (d+6)}{6\alpha-d}}$, $\alpha=\frac{d(d+2)}{d^2-d-2}\ge 2$,
$\beta=\frac{-d^2+4d+4}{4d}$, $d\alpha/(d+\alpha)=1+2/d,$ to get the 
estimate by the time integral of $|n|_{1+2/d}$ but this time for small mass only. 
To get a bound for $|n|_{1+2/d}$, assumptions have to be
more restrictive, e.g. Lemma \ref{le4p} requires $2\le d\le 4$.   

In fact, we have obtained the estimate
\be \label{es3pp}
\vert n(t)\vert_2^2 + c\int_0^t \theta P'^2 |\nabla n(s)|_{2}^2\ ds \le C\,,
\ee
which due to the estimate $P'\ge p_0$ implies \rf{es3}. 
\qed

\medskip

\noindent{\bf Proof of the global existence part of Proposition \ref{loc}}.
We are now ready to prove (\ref{global}) and thus obtain the global
existence. Let $T>0$ and $t\in [0,T_{max})\cap [0,T]$. We
claim that there is $C>0$ depending on
$n_0$ and $T$ and bounds on $\th$ such that  
\be
\label{tt3}
\vert nP'\vert_{L^{2(d+4)/(d+2)}(\Omega\times (0,t))} + \vert\nabla
(nP')\vert_{L^2(\Omega\times (0,t))}\le C\,. 
\ee
Indeed, we infer from assumption \rf{p0}-\rf{p3} that
$zP'(z)\le C\ (1+z^{1+2/d})$ and $zP''(z)\le C\ (1+z^{2/d})$ for $z\ge 0.$ 
Consequently,
\be\label{PH1}
\sup_{s\in [0,t]} \vert n(s)P'(s)\vert_{\frac{2d}{d+2}} 
+ \int_0^t \Vert n(s)P'(s)\Vert_{H^1}^2\ ds \le C\, 
\ee
holds. Next, we use the continuity of the imbedding of $H^1(\Omega)$ in
$L^{2d/(d-2)}(\Omega)$ and an interpolation argument to deduce (\ref{tt3}). 

We now employ a bootstrap argument to show that (\ref{global}) holds
true. It follows from (\ref{p1}) and the Sobolev imbedding that 
$$
\vert n\vert^{(d+2)/d}_{L^{2+4/d}(0,t;L^{2(d+2)/(d-2)}(\Omega))} 
\leq C\vert \th^{d/2+1}P\vert_{L^2(0,t;L^{2d/(d-2)}(\Omega))} \le C\,,
$$
which, together with (\ref{es3}), leads to
$$
\int_0^t \vert n(s)\vert_{2+8/d}^{2+8/d}\ ds \le \int_0^t \vert
n(s)\vert_{2(d+2)/(d-2)}^{2+4/d}\ \vert n(s)\vert_2^{4/d}\ ds \le C\,.
$$
Therefore,
$$
\vert n\vert_{L^{2+8/d}(\Omega\times (0,t))} \le C\,,
$$
and we infer from \cite[Theorem~IV.9.1 and
Lemma~II.3.3]{LSU} that  
$$
\vert\nabla\varphi\vert_{L^{2(d+2)(d+4)/(d^2-8)}(\Omega\times (0,t))} +
\vert\Delta\varphi\vert_{L^{2+8/d}(\Omega\times (0,t))} \le C\,.
$$
This estimate and (\ref{tt3}) ensure that
$$
\vert \nabla (nP')\cdot \nabla\varphi\vert_{L^{(d+4)(d+2)/(d^2+3d)}(\Omega\times
(0,t))}
+ \vert nP'\ \Delta\varphi\vert_{L^{(d+4)/(d+2)}(\Omega\times (0,t))} \le
C\,. 
$$
Since
\be\label{nFe}
n_t - \th\nabla \left(P'^2 \nabla n\right)= \nabla (nP')\cdot \nabla\varphi + nP'\
\Delta\varphi\,, 
\ee
we use once more \cite[Theorem~IV.9.1]{LSU} to obtain that 
$$
\Vert n\Vert_{W_{(d+4)/(d+1)}^{2,1}(\Omega\times (0,t))} \le C\,,
$$
which, in turn, implies that 
$n\in L^{\frac{(d+2)(d+4)}{d^2+d-6}}(\Omega\times (0,t)).$ 
With thus improved $n$ we would like to bootstrap once again.
The  right hand side of (\ref{nFe}) is in the space 
$L^{q}(\Omega\times (0,t))$, with 
$$q=d(d+2)(d+4)\min \left\{1/(d^3+4d^2-12),1/(2(d+1)(d+3)(d-2)) \right\}\,.$$ 
Therefore, for $d=3$ we finally get
the right hand side of \rf{nFe} in $L^{q}(\Omega\times (0,t))$ with $q=35/12$ larger 
than critical $1+3/2$ allowing to conclude with
\be \label{hol}
\Vert n\Vert_{\mathcal{C}^\varepsilon([0,t])} \le C\,,
\ee
for some $\varepsilon>0$ by \cite[Lemma~II.3.3]{LSU}. 
However, for $d=4$ one should note that 
we have obatined from the boostrap the integrability 
of the  right hand side of (\ref{nFe}) of the order $q=48/35$
which is less than we had before, i.e. $q=8/5$. Thus for $d=4$
we cannot conclude with the estimate \rf{hol}. \qed

\section{The local existence result for the original elliptic--parabolic problem}

In this section we shall subtract a convergent subsequence of solutions 
to \rf{reg1}--\rf{reg4} obtained in the previous section which will guarantee
the following existence result for the limiting problem \rf{ibp1}--\rf{ibp4}
as $k\rightarrow\infty$.

\begin{thm}\label{thglobal}
Assume that $M$ is small enough if necessary, 
and $P$ is smooth and satisfies \rf{p0}--\rf{p00}
and \rf{PRp}. Moreover, let
$2\le d \le 4$, $n_0\in L^{2}(\Omega)$ and, for given 
constants $0<a<b$, $\th\in C(0,T;[a,b])$. 
Then there exist a weak local-in-time solution 
$n\in {\cal C}\left(0,T; 
L_w^2(\Omega)\right),$ $\f\in L^\infty\left(0,T; 
H^2(\Omega)\right),$ $\th^{d/2+1}P\in
L^2\left(0,T;H^1(\Omega)\right)$
of the system \rf{bvp1}--\rf{bvp4}, i.e.
\bea
\iom{(n-n_0)\, \chi} &+& 
\int_0^t \iom {\nabla\chi\cdot\left(\th P'^2\nabla n 
+nP'\nabla\f\right)}ds =0,\label{xy}\\
\Delta\f&=&n, \ \ \ \f=0\ \ {\rm {on}}\ \partial\Omega, \label{xz}
\eea
for each test function $\chi\in W^{1,2d/(d-2)}(\Omega)$.
Additionally,
\bea
\label{yx}
& & \vert n(t)\vert_{1+2/d} + \Vert\varphi(t)\Vert_{H^1} \le C\,,\\
\label{yz}
& & \vert n(t)\vert_2 + \int_t^{t+1} \Vert 
\th^{d/2+1}(s)P(s)\Vert_{H^1}^2\ ds \le C\,,
\eea
for any $t\in [0,T)$ where $C$ depends on $n_0$, $\Omega$, 
$d$, $a$ and $b$. If $d\le 3$ then the global result can be claimed.
\end{thm}
\noindent\textbf{Proof.}
We follow the lines of the proof from \cite{BLN}, where $d=3$ and a constant 
temperature $\th$ were assumed. We consider $n_0\in L^{2}(\Omega)$ such $n_0\ge 0$ 
a.e. in
$\Omega$ and put $M=\vert n_0\vert_1$ (sufficiently small if necessary). 
Let  $(n_{0,k})_{k\ge 1}$ be a~sequence of nonnegative functions in
$\mathcal{C}^\infty(\bar{\Omega})$ approximating $n_0,$ i.e., 
\be
\label{cvinit}
\vert n_{0,k}\vert_1=M \;\;\mbox{ and }\;\; \lim_{k\to \infty} \vert
n_{0,k} - n_0\vert_{2} = 0\,.
\ee
For $k\ge 1$, we denote by $(n_k,\varphi_k)$ the unique classical
solution to (\ref{reg1})--(\ref{reg4}) with initial datum
$(n_{0,k},0)$ given by Theorem~\ref{loc} and let $P_k=P(n_k\th^{-d/2})$. 
Owing to \rf{le4p}, (\ref{es3pp}) and (\ref{cvinit})
there is $C>0$ such that 
\be\label{u}
 \vert n_k\vert_{2} + \Vert\varphi_k\Vert_{H^1} +
\frac{1}{k}\ \int_0^T \vert(\varphi_{k})_s(s)\vert_2^2 \ ds
+ \int_0^T \|\th^{d/2+1}(s)P_k(s)\|_{H^1}^2\, ds  \le C\,.
\ee
Observe that the H\"older inequality, \rf{u} and assumptions \rf{p0}, (\ref{p2}) 
imply 
\bean
&  \vert\th P_k'^2 \nabla n_k\vert_{\frac{2d}{d+2}} &\le C\left(\vert n_k^{2/d}
\ \nabla
P_k\vert_{\frac{2d}{d+2}}+\th^{d/2+1}\vert \nabla
P_k\vert_{\frac{2d}{d+2}}\right) \\
& & \le C\left( \vert n_k \vert_{2}
\vert \nabla P_k\vert_{2}+\vert \nabla
P_k\vert_{\frac{2d}{d+2}}\right),\\
&  \vert n_kP'_k\nabla \f_k \vert_{d/2} &\le C  \vert P_k\th^{d/2}\
\nabla\varphi_k\vert_{d/2} \le C  \vert P_k\vert_{\frac{2d}{d-2}}
\vert\nabla\varphi_k\vert_{\frac{2d}{6-d}},
\eean
whence, by \rf{u} thanks to the imbedding of $H^1(\Omega)$ in
$L^{\frac{2d}{d-2}}(\Omega)$,
\be
\label{tt5}
\int_0^T \left( \vert\th(s) P'^2_{k}(s)\nabla n_k(s)\vert_{2d/(d+2)}^2 + \vert 
n_k(s)P'_k(s)\ \nabla\varphi_k(s)\vert_{d/2}^2 \right)\ ds \le C\,.
\ee
For $d=2$ due to (\ref{fii}) we get the $L^{\infty}$ bound for $\nabla 
\varphi_k$ whence $|n_kP'_k\nabla\varphi_k|_1\le c|P_k|_1\le C$.
We then deduce from the above inequality and equation \rf{reg1} that
\be
\label{tt6}
\vert (n_k)_t\vert_{L^2(0,T;W^{1,2d/(d-2)}(\Omega)')} \le C\,.
\ee
Consequently, owing to (\ref{tt6}), \rf{es3} and \rf{u} the sequence $(n_k)$ is bounded in 
$L^2(0,T;H^1(\Omega))$ and in $H^1(0,T;W^{1,2d/(d-2)}(\Omega)')$. Owing to the 
compactness of the imbedding of
$H^1(\Omega)$ in $L^2(\Omega)$ and to the continuity of the imbedding
of $L^2(\Omega)$ in $W^{1,2d/(d-2)}(\Omega)'$, we infer from
\cite[Corollary 4]{S} that $(n_k)$ is relatively compact in
$L^2(\Omega\times (0,T))$. Therefore, there are $n\in
L^2(\Omega\times (0,T))$ and a subsequence of $(n_k)$ 
such that $n_k\rightarrow n$ a.e. and
\be
\label{cvn}
n_k \longrightarrow n \;\mbox{ in }\; L^2(\Omega\times
(0,T))\cap \mathcal{C}([0,T];W^{1,2d/(d-2)}(\Omega)')\,.
\ee
 Let $\varphi\in L^\infty(0,T;H^2(\Omega))$ be the solution to
\be
\label{eqphi}
\Delta\varphi = n \;\;\mbox{ in }\;\; \Omega\times (0,T)\,, \qquad
\varphi=0 \;\;\mbox{ on }\;\; \partial\Omega\times (0,T)\,.
\ee
It follows from (\ref{reg2}) and (\ref{eqphi}) that $\varphi_k -
\varphi$ solves the Poisson equation
$$
-\Delta(\varphi_k - \varphi) = n-n_k - \frac{1}{k}\ (\varphi_k)_t
$$
with the homogeneous Dirichlet boundary conditions, and the right-hand
side of the above equation converges to zero in $L^2(\Omega\times
(0,T))$ by (\ref{u}) and (\ref{cvn}).
Therefore,
\be
\label{cvphi}
\varphi_k \longrightarrow \varphi \;\;\mbox{ in }\;\;
L^2(0,T;H^2(\Omega))\,.
\ee
Combining (\ref{u}) with the convergence results (\ref{cvn}) and
(\ref{cvphi}) finally allow us to conclude that $P_k'^2\nabla n_k$ and
$P'_k n_k \nabla\varphi_k$ converge weakly to $P'^2 \nabla\ n$ and
$nP'\nabla\varphi$ in $L^{2d/(d+2)}(\Omega\times (0,T))$ and
$L^{d/2}(\Omega\times (0,T))$, respectively.
It is now straightforward to pass to the limit as $k\to \infty$ 
and conclude that $(n,\varphi)$ is a
weak solution to (\ref{bvp1})--(\ref{bvp3}) as stated in
Theorem \ref{thglobal}.

We may also pass to the limit in (\ref{yx})
and use classical lower semicontinuity argument to deduce that
(\ref{yx}) holds true.

Next, by \rf{p0} and \rf{p1} it follows from the conservation of mass,
\rf{dnl} and the Poincar\'e inequality that
\be \label{dnk}
\frac{d}{dt} \vert n_k(t)\vert_2^2 + \gamma\ \left( \vert n_k\vert_2^2 +
  \Vert \th^{d/2+1}P_k\Vert_{H^1}^2 \right) \le C\ \left( 1 +
  \vert\nabla\varphi_k\vert_{\alpha}^{\frac{2\alpha(d+6)}{6\alpha-d}} \right)
\ee
for some positive constant $\gamma$. 
Integrating with respect to time, we get
\be
\label{lila}
\vert n_k(t)\vert_2^2 \le \vert n_{0,k}\vert_2^2\ e^{-\gamma t} + C\
\int_0^t \left( 1 + 
\vert\nabla\varphi_k(s)\vert_{\alpha}^{\frac{2\alpha(d+6)}{6\alpha-d}} \right)\
e^{\gamma (s-t)}\ ds
\ee
for $t\ge 0.$ Next, from the Fubini theorem and the double integration of (\ref{dnk}) 
we obtain, for $t\ge 1,$
\bea
& & \int_t^{t+1} \Vert \th^{d/2+1}(s)P_k(s)\Vert_{H^1}^2\ ds \le
\int_{t-1}^t \int_\tau^{\tau+2} \Vert \th^{d/2+1}(s)P_k(s)\Vert_{H^1}^2\
ds \, d\tau  \nonumber\\
\label{lili}
& & \le C\ \left( 1 + \int_{t-1}^t \vert
  n_k(\tau)\vert_2^2\ d\tau + \int_{t-1}^{t+2}
  \vert\nabla\varphi_k(s)\vert_{\alpha}^{\frac{2\alpha(d+6)}{6\alpha-d}}\ ds \right).
\eea
Now, $\vert\nabla\varphi_k\vert_{\alpha}$ is bounded in
$L^q(0,t)$ for any $q\in (1,\infty)$ by Lemma~\ref{le2}, and we
infer from (\ref{cvphi}) and the continuous imbedding of
$H^2(\Omega)$ in $W^{1,\alpha}(\Omega)$ that $\vert \nabla\varphi_k
-\nabla\varphi\vert_{\alpha}$ converges to zero in
 $L^2(0,t)$. Consequently, by interpolation, $\vert \nabla\varphi_k
  -\nabla\varphi\vert_{\alpha}$ converges to zero in $L^{\frac{2\alpha(d+6)}{6\alpha-d}}(0,t)$. 
Then one can pass to the limit as $k\to \infty$ in (\ref{lila}) and
  (\ref{lili}) with the help of (\ref{cvn}) and weak convergence
  arguments for the left-hand sides and conclude that
$$
\vert n(t)\vert_2^2 \le \vert n_0\vert_2^2\ e^{-\gamma t} + C\
\int_0^t \left( 1 + 
\vert\nabla\varphi(s)\vert_{\alpha}^{\frac{2\alpha(d+6)}{6\alpha-d}} \right)\
e^{\gamma (s-t)}\ ds
$$
for $t\ge 0,$ while for $t\ge 1$
$$
\int_t^{t+1} \Vert \th(s)^{\frac{d+2}{2}}P(s)\Vert_{H^1}^2 ds 
\le C \left( 1 + \int_{t-1}^t \vert n(s)\vert_2^2 ds + \int_{t-1}^{t+2} 
  \vert\nabla\varphi(s)\vert_{\alpha}^{\frac{2\alpha(d+6)}{6\alpha-d}}ds \right)\,.
$$
Since $\varphi$ is a solution to equation
(\ref{bvp2}), by (\ref{yx}) we get taking $\alpha=\frac{d(d+2)}{d^2-d-2}$ 
$$
\vert\nabla\varphi\vert_{\alpha} \le C\
\Vert\varphi\Vert_{W^{2,1+2/d}} \le C\ \vert n\vert_{1+2/d}\le C.
$$
Inserting this estimate in the previous two inequalities yields the 
boundedness of $\vert n(t)\vert_2$ with respect to time,
and then (\ref{yz}). For $d=2$ we have $1+2/d=2$ whence the estimate
for  $\vert n\vert_{1+2/d}$ is sufficient.\qed

\section{Fixed point for the temperature operator ${\cal{T}}$}

First we recall a lemma on relations between $n$ and $\vt$ imposed by \rf{ene}.
This should be understood as necessary condition for the
density obtained from \rf{bvp1}-\rf{bvp4} and not as a~sufficient condition
for admissibility of the given energy $E$. The lemma on {\emph{a~priori}} bounds is 
related to the one from \cite{BLN} in 
the Fermi--Dirac case and to the ones from \cite{BS-d} and \cite{BS-t} in more general case. 
Recall from \cite[Lemma 3.1]{BS-t} or \cite{S-b} the following version of these 
energy estimates.
\begin{lem}\label{Ebd}
Let $\nu=4/(d(4-d))$ for any $2\le d<4$. Provided that $P(s)\ge p_1 s^{1+2/d}$ for any
$\frac{d}{2}p_1>\varepsilon>0$ and all $s\geq 0$, 
the following estimate holds
\begin{equation}\label{nes}
E+CM^{1+\nu} \ge 
\max\left\{\varepsilon\iom{n^{1+2/d}},\left|\f\right|_2^2\right\}.
\end{equation}
Moreover, for each $0<\varepsilon<d/2,$ the temperature $\vt$ and the density $n$ should satisfy 
\begin{equation}\label{pes}
E\ge \varepsilon \iom{\vt^{d/2+1}P(n\vt^{-d/2})}+\left|\iom{n\f}\right|-CM^{1+\nu}.
\end{equation}
\end{lem}

Now we shall prove some {\emph{a~priori}} estimate for $L^{1+2/d}$ norm
of the solution to BVP (\ref{bvp1})--(\ref{bvp4}). We derive them directly
from these equations since at this moment we cannot directly use the energy
{\emph{a~priori}} bounds presented above. Note that these 
{\emph{a~priori}} estimates for limit functions are better  
than those for the perturbed parabolic system presented
in previous sections (cf. lemmas: \ref{le3p}, \ref{Pbd}, \ref{le4p}).

The next lemma can be found in \cite{S-b} (cf. Lemma 2.1 and 2.2 therein).

\begin{lem}\label{le3} 
For any $2\le d< 2(1+\sqrt{2})$ and $\varphi$ related to $n$ by \rf{bvp2} we have the estimate
\begin{equation}\label{fie}
\left|\iom{n\f}\right|\le
CM^{1-2/d}\iom{n^{1+2/d}}.
\end{equation}
Let $2\le d\le 4$ and assume that 
\be \label{PR}
P(z)=p_1z^{1+2/d}+R(z),
\ee
where the lower order term satisfies $\vert R'(z)\vert z^{1/2-1/d}\leq B.$ 
Define the `asymptotic energy', i.e. $E^a(t)=\lim_{\vt\r 0^+} E(t)$, by 
\be \label{Eas}
E^a(t)=\frac{d}{2} \iom{p_1\,n^{1+2/d}} - \frac{1}{2} \iom{\vert \nabla \varphi\vert^2}.
\ee
Then for any fixed $T>0$ and any $t\in [0,T]$ the following 
growth condition for $E^a$ is available
\be
\label{es22}
\frac{d}{dt}E^a(t) \leq C\theta (t)^{d/2}  \iom{\vert \nabla \varphi\vert^2}.
\ee
\end{lem}

\medskip

\noindent {\bf Remark.} Note that the above theorem holds both 
in the polytropic case with $R(z)=0$ and in the Fermi--Dirac case,
since by the properties of Fermi function (cf. Lemma \ref{FAS} the Appendix and 
\cite[Sec.5]{BNS}) $\vert R'(z)\vert \leq Bz^{-2/d}$
at $z=\infty$ and $\vert R'(z)\vert \leq B$ at $z=0$.

\medskip

\noindent {\bf Remark.} It should be noted that for the polytropic case the theorem 
implies the dissipation of the energy, since in this case $E^a=E$.

\medskip 

Applying the estimate (\ref{fie}) to $\iom{n \varphi}= -\iom{\vert \nabla \varphi\vert^2}$
and integrating (\ref{es22}) from Lemma 
\ref{le3} allows us to derive the following corollary (for details and the proof 
see \cite{S-b}).

\begin{cor}\label{co1}
Under the assumptions of Lemma \ref{le3} $E^a$ grows like
\be \label{Eag}
E^a(t)\leq E^a(0) \,\exp{(C(t))}
\ee
where the function $C$ is defined by 
$C(t)=cM^{1-2/d}\int_0^t(\th(s))^{d/2}ds$. Moreover, $E^a$ is positive
if $CM^{1-2/d} <dp_1$ while for $d=2$ we assume smallness of the 
Poincar\'e constant i.e. $C<2p_1$.
\end{cor}

After integrating inequality (\ref{Eag}) from Corollary \ref{co1}
and using Lemma \ref{Ebd} we obtain $L^{1+2/d}$ estimate for the density $n.$

\begin{cor}\label{co2}
Under the assumptions of Lemma \ref{le3} we have for any $2\le d<4$, $\nu=4/(d(4-d))$ and any $M>0$
$$
\iom{n^{1+2/d}} \leq CM^{1+\mu}+\vert E^a(0)\vert \exp{(C(t))},
$$
while for small $M>0$ and any $2\le d\le 4$ the constant $C$ may depend on $M$
$$
\iom{n^{1+2/d}} \leq C E^a(0) \exp{(C(t))}
$$
where $CM^{1-2/d}<dp_1$.

\end{cor}

\medskip

Corollaries \ref{co1} and \ref{co2} allow us to 
to define the new 
temperature $\vt$. This was the subject of the considerations in \cite{S-b} 
under physically acceptable property of the pressure $\frac{\partial p}{\partial 
\vt}>0$ expressed as
\be \label{pse}
P(z)z^{-1-2/d}\searrow p_2>0
\ee
that guarantees, in particular, 
the uniqueness of the temperature $\vt$ emerging from the energy formula (\ref{E}). 

The next theorem claims that that temperature is well defined for some values of the energy 
(for the proof see Theorem 3.2 in \cite{S-b}),
and the remainder of the section is devoted to proving its compactness.

\begin{thm}
Assume that $P$ is smooth and satisfies (\ref{p0}), (\ref{p00}), \ref{pse}) and \rf{PRp}
Then the temperature operator
$\mathcal{T}:\th\mapsto\vt$ is formally well defined by \rf{E} for small mass 
$M\ll 1$ and $2\le d \le 4$ all the values of the energy $E$ admissible at $t=0$. Moreover, 
for $2\le d<4$, $\nu=4/(d(4-d))$ and some positive constants $B,C$, it has to satisfy
\be\label{Efi}
E<BM^{1+2/d}-CM^{1+\nu}.
\ee
\end{thm}

Next, we estimate $\vt'$ to get the compactness of the operator ${\mathcal T}$.
By differentiation of the energy relation \rf{E} we get
\be \label{vtd}
\vt'=\left(\iom{\frac{\partial p}{\partial n}n_t}+\frac{1}{d}\frac{d}{dt}
\iom{n\f}  \right)\left( \iom{\frac{\partial p}{\partial \vt}} \right)^{-1}.
\ee
In the following two lemmas we claim the boundedness of both factors
in appropriate norms so that $\vt'$ be in $L^\gamma$ with some
$\gamma>1$ which guarantees the equicontinuity condition in
the classical Arz\`ela-Ascoli thoerem.

\begin{lem}\label{inp}
Assume that $P$ is a smooth function such that 
\bea
\left(P(z)z^{-1-2/d}\right)'<0\label{Pp},\\
\left((P(z)z^{-1-2/d})'z^{1+2/d}\right)'> 0.\label{Ppp}
\eea
Then the first inequality implies \rf{p2} with a strict inequality and $C=1+2/d$, i.e.,
\bea \label{Pd12}
P'(z)z< (1+2/d)P(z)\,.
\eea
Moreover, the function $p$ is decreasing, convex with respect to 
$\vt$ and satisfies
\be\label{pfb}
\iom{\frac{\partial p}{\partial \vt}}>C
\ee
for some $C>0$ depending on $M$ and a lower bound for $\vt\ge a$
provided that
\bea \label{Pd34}
\liminf_{z\rightarrow\infty}\left(-P(z)z^{-1-2/d}\right)'z^{4/d+1} \ge C>0,\\
\liminf_{z\rightarrow 0}\left(-P(z)z^{-1-2/d}\right)'z^{2/d+1} \ge C>0.\label{Pd35}
\eea
 \end{lem}
\noindent{\bf Proof.} The formula for the first derivative reads
\be \label{Pd1}
\frac{\partial p}{\partial \vt} = -(d/2)\vt^{d/2}\left( P'(z)z-(1+2/d)P(z) \right),
\ee
or in an another form
\be
 \frac{\partial p}{\partial \vt} = -(d/2)\vt^{d/2}z^{2+2/d}
\left( P(z)z^{-1-2/d}\right)'>0\,,
\ee
where $z=n\vt^{-d/2}.$ Then, the second derivative can be calculated 
\be\label{Pd2}
(2/d)^2\vt^{1-d/2} \frac{\partial^2 p}{\partial \vt^2}
= P''(z)z^2-(1+2/d)P'(z)z+(1+2/d)P(z)\,,
\ee
or expressing it in a more concise way
\bean
(2/d)^2\vt^{1-d/2} \frac{\partial^2 p}{\partial \vt^2} 
=  z^2 \left((P(z)z^{-1-2/d})'z^{1+2/d}\right)'\,. 
\eean
Thus the convexity of $p$ with respect to $\vt$ follows from the second 
assumption (which by the way can be deduced from the first assumption or
(\ref{Pd12}) under an extra convexity assumption on $P$).
Now, by the asymptotics of $P$, i.e. (\ref{Pd34}), \rf{Pd35},  
$\liminf_{n\rightarrow 0}\iom{\frac{\partial p}{\partial \vt}}\geq CM$ as and 
$\liminf_{n\rightarrow\infty}\iom{\frac{\partial p}{\partial \vt}} \geq CM^{1-2/d}a,$
respectively. Hence, by the convexity of $p$ with respect to $\vt$, \rf{pfb} follows. 
\qed

\medskip

To get the bound for $\vt'$ we are left to estimate the denominator in (\ref{vtd}).

\begin{lem}\label{vtn}
The denominator appearing in \rf{vtd} is bounded in some $L^{\gamma}$ with $\gamma>1$, i.e.,
$$\int_0^T \left|\iom{\frac{\partial p}{\partial n}n_t}+\frac{1}{d}\frac{d}{dt}
\iom{n\f}\right|^\gamma\,dt<C\,.$$
\end{lem}
\proof
First, recall that 
$$
\iom{n_t\f }=\frac{d}{dt}\left( \iom{n\f +\frac{1}{2}\vert \nabla \varphi\vert^2} \right)
+ \frac{1}{k}\iom{\varphi_t^2}.
$$
Using (\ref{tt6}) we get the bound for $n_t$ in $L^2(0,T;W^{1,2d/(d-2)}(\Omega)')$ so
we have to show that both $\frac{\partial p}{\partial n}$ and $\f$ are bounded in 
$L^\xi(0,T;W^{1,2d/(d-2)}(\Omega))$ with some $\xi>2$. Using (\ref{le2}) and \rf{yz}
with $\alpha=\frac{2d}{d-2}$ and any $q>1$ 
we get 
$$\int_0^t |\nabla\f(s)|_{\frac{2d}{d-2}}^q\, ds \le C \int_0^t |n(s)|_2^q ds\, \le C'\,.$$
whence $\f\in L^\xi(0,T;W^{1,2d/(d-2)}(\Omega))$ with $\xi>2$.
Moreover, $\frac{\partial p}{\partial n}=\th P'\sim n^{2/d}$ and the function $\nabla \frac{\partial p}{\partial n}=
\th^{1-d/2}P''$ is bounded in view of the regularity assumption on $P$, 
so the claim is guaranteed by the estimates for $P$ 
\rf{p0} and 
\rf{p00}
and \rf{le2}.
Lastly, 
from (\ref{es3p}) follows the bound for $\iom{\varphi_t^2}$ and from 
(\ref{pfi}) and Lemma \ref{le4p} for $\iom{\vert \nabla \varphi\vert^2}$. 
\qed

Finally, we recall after \cite{BS-t} and \cite{S-b} 
{\emph{a~priori}} bounds on the fixed points of the compact operator $\mathcal{T}$ thus 
guaranteeing the existence result for the problem.

The authors assumed therein, for negative initial values of the entropy
\be\label{Wen}
\liminf_{z\rightarrow\infty} \left( H(z)-(d/2+1)P(z)/z \right) > {\mathcal W}(0)/M\, ,
\ee
where $H'(z)z=P'(z)$ under the following conditions consistent with \rf{p0}
\bea
P(z)/z^{1+2/d}\searrow \varepsilon>0,\label{Pz1}\\
\liminf_{z\searrow 0}P(z)/z>0.\label{Pz2}
\eea
If \rf{PRp} is fulfilled then the highest order
terms cancel assuming that the limit exists
\bea\label{G0}
\lim_{z\rightarrow\infty} \frac{ H(z)z-(d/2+1)P(z)}{z}
= \lim_{z\rightarrow\infty} \left( H(z)-\frac{d}{2} zH'(z) \right) \\
= \lim_{z\rightarrow\infty} \left( G(z)-(d/2) zG'(z) \right)\stackrel{\rm df}=G_0
\eea
where
\be
H'(z)=h_1(2/d) z^{2/d-1} + G'(z)
\ee
with $G'(z)=o(z^{2/d-1})$, we are left with the analysis of the lower 
order term $G'(z)=g_1z^{\beta}+o(z^{\beta})$ with some $\beta < 2/d-1$.
Namely, if $\beta<-1$ then $G_0=0$ in \rf{G0}, e.g. for the Fermi--Dirac 
model in 
$d=2$. Otherwise, if $\beta \in [-1,2/d-1)$ then the important factor is the sign of
$g_1$ which has to be positive, and indeed is, e.g. $g_1=1-2/d$ 
and $\beta=-2/d$ in the Fermi--Dirac 
case for $d\geq 3$, to imply $G_0=\infty$ and thus to guarantee \rf{Wen}.

Thus we have proved the existence of a fixed point for temperature operator $\mathcal{T}$
and we can formulate the following existence result in the microcanonical case.

\begin{thm}\label{final}
Assume that $P$ is smooth, satisfies \rf{p0}, \rf{p00}, \rf{PRp}, \rf{pse}, \rf{Pz1}, \rf{Pz2} and
\rf{G0} with $G_0\ge 0$.
Then for the negative initial values of the entropy we get the global existence result 
for \rf{bvp1}-\rf{bvp4} with the energy constraint \rf{E}.
\end{thm}

\section{Appendix on Fermi--Dirac model}

First, it should be noted that for the Fermi--Dirac case we have
\be \label{PFD}
dP_{FD}(z)=\mu f_{d/2}\left(f_{d/2-1}^{-1}(2z/\mu)\right),
\ee
where $\mu =\eta_0 G \sigma_d 2^{d/2}$, $G$ is the gravitational constant, $\eta_0$ --
a bound for the density in phase space and $f_{\alpha}$ is the Fermi function of order
$\alpha>-1$ defined by 
\be \label{Ff}
f_{\alpha}(z)=\int_0^\infty\frac{x^{\alpha}}{1+e^{x-z}}dx.
\ee
 
\medskip

In \cite[Lemma 5.1]{BNS} substitute $z=-\log(\lambda)$
and $f_{\alpha}(z)=I_{\alpha}(e^{-z})$ to get

\begin{lem}\label{exas}
The following asymptotic relations hold as $z\rightarrow \infty$
\be
f_\alpha(z)-\frac{z^{\alpha+1}}{\alpha+1}
={\cal O}\left(z^{\alpha-1}\right) ,\label{r}
\ee
for each $\alpha\ge 0$, while for each $\alpha>-1$
\be
z^{-\alpha}\left\{f_\alpha(z)z
-\frac{\alpha+2}{\alpha+1}f_{\alpha+1}(z)\right\}\to-\frac{\pi^2}{3}\,.
\label{R}
\ee
Moreover, we have the recursive relation for the derivatives
\be\label{i}
f'_\alpha(z)=\alpha f_{\alpha-1}(z).
\ee
\end{lem}

\medskip

Next \cite[Lemma 2.2]{BLN} can be formulated as follows.

\begin{lem}\label{fab}
For $\alpha>\beta,$ $f_\alpha \circ f_\beta^{-1}$ is an increasing convex function.
\end{lem}

In conclusion of the above lemma the function $P_{FD}$ shares the same properties.

\begin{lem}\label{Pic}
The function $P_{FD}$, defined by (\ref{PFD}), is increasing and convex function.
\end{lem}

Next, to check that the assumptions of the Lemma \ref{inp} are verified
for the Fermi--Dirac, case we will need a version of \cite[Lemma 5.3]{BNS}.

\begin{lem}
For all $\beta<\alpha+1$ the following inequality holds
\be\label{fin}
f'_{\alpha+\beta}f'_{\alpha-\beta}-(f'_\alpha)^2>0\,.
\ee
\end{lem}

\medskip

Now, we recall after \cite[Lemma 2.1]{BLN} the following properties of the 
Fermi--Dirac pressure.

\begin{lem}\label{DVF}
The function $P_{FD}$ belongs to ${\cal C}^2([0,\infty))$, is nonnegative,
increasing, convex and can be extended to an element (still denoted by $P_{FD}$)
of ${\cal C}^2([-\delta,\infty))$ for some $\delta>0$.
\end{lem}

Next the following asymptotic result holds (cf. \cite{S-b}) for $$P_{FD}(z)=p_1 z^{1+2/d}+R_{FD}(z).$$

\begin{lem}\label{FAS}
For the Fermi--Dirac pressure $P_{FD}$ we have at $z=\infty$
\be\label{FDi}
P_{FD}'(z)=\frac{2}{d}\left(\frac{d}{\mu}z\right)^{2/d}+{\mathcal O}(z^{-2/d})
\ee
and, in consequence,
\be
P_{FD}(z)=p_1 z^{1+2/d}+{\mathcal O}(z^{1-2/d}).
\ee
where $p_1=\frac{2}{d+2}(d/\mu)^{2/d}.$ Moreover, at $z=0$, we have
\be \label{FD0}
R_{FD}'(z)={\mathcal O}(1)\,.
\ee
\end{lem}

\medskip

\begin{lem}
The conditions $\frac{\partial p}{\partial \vt}>0$ and/or $P(z)/z^{1+2/d}\searrow 
p_2>0$ are satisfied for 
\be
P_{FD}(z)=(\mu/d) f_{d/2}\circ f_{d/2-1}^{-1}\left(2z/\mu\right).
\ee
\end{lem}
\proof
Indeed, putting $2z=\mu f_{d/2-1}(x)$, we get
$$
\frac{\partial p}{\partial \vt}=
\frac{d}{2}\vt^{d/2}\left((-\mu/2)(f_{d/2-1}(x))^2(f'_{d/2-1}(x))^{-1}
+(1+2/d)(\mu/d)f_{d/2}(x)\right)\,.
$$
Thus, the condition $\frac{\partial p}{\partial \vt}>0$ is equivalent to
$$
-\left((\mu/2)f_{d/2-1}(x)\right)^2+(\mu/d)f'_{d/2-1}(x)(\mu/2)(1+2/d)f_{d/2}(x)>0\,.
$$
This, however, follows from Lemma \ref{fab} (take $\alpha=d/2$ and $\beta=d/2-1$) 
or more explicitly
by the property of Fermi functions presented in \rf{fin}, namely
$
\frac{d}{dx}\left(\frac{f_{d/2}(x)}{f_{d/2-1}(x)}\right)<0\,. 
$
\qed

Lastly, we shall trace how (\ref{bvp1})--(\ref{bvp4}) could 
be derived, in the Fermi--Dirac case, from (\ref{ibp1})--(\ref{ibp4}) 
(used by the authors in \cite{BLN}) under the assumption
(\ref{ssf}) with a specific diffusion coefficient, 
used in \cite{BLN} and \cite{CRS},
\be\label{D}
D(\lambda)=\frac{-I_{d/2-1}(\lambda)}{\lambda I'_{d/2-1}(\lambda)},
\ee
where $I_\alpha(e^{-z})=f_\alpha(z)$ and 
$\lambda = I_{d/2-1}^{-1}\left(\frac{2n}{\mu\th^{d/2}}\right).$
Using the recurrence property \rf{i} (cf. also \cite[Section 5]{BNS} and
\cite[Lemma 1.1]{S-d}), we get, differentiating formula \rf{PFD},
the relation $D=P'$. Furthermore, $ \th D= \frac{\partial p}{\partial n}.$
Moreover, it should be noted that in \cite{BLN} and \cite{BNS}
the authors used the following notation $F'=P'^2, V=nP'.$ Note that $D$ should
be defined exactly as in (\ref{D}) but it might differ throughout these papers up to
a constant, inessential therein.


\end{document}